\documentclass[reqno]{amsart}
\usepackage{amsfonts,amssymb}
\usepackage{enumerate}

\newtheorem{thm}{Theorem}
\newtheorem{lem}[thm]{Lemma}

\theoremstyle{definition}

\theoremstyle{remark}

\author{Jing-Jing Huang \and Robert C. Vaughan}
\address{
JJH: Department of Mathematics, McAllister Building,
Pennsylvania State University, University Park, PA 16802-6401,
U.S.A.}
\email{huang@math.psu.edu}
\address{
RCV: Department of Mathematics, McAllister Building,
Pennsylvania State University, University Park, PA 16802-6401,
U.S.A.}
\email{rvaughan@math.psu.edu}
\thanks{JJH: Part of this paper is prepared when the first author was visiting the National Center for Theoretical Sciences in
Taiwan. Their hospitality and financial support are gratefully acknowledged by him.}  
\thanks{RCV: Research supported in part by NSA grant number H98230-09-1-0015.}

\begin{document}

\title
{Mean Value Theorems for Binary Egyptian Fractions II}

\begin{abstract}
In this article, we continue with our investigation of the Diophantine equation $\frac{a}n=\frac1x+\frac1y$ and in particular its number of solutions $R(n;a)$ for fixed $a$. We prove a couple of mean value theorems for the second moment $(R(n;a))^2$ and from which we deduce $\log R(n;a)$ satisfies a certain Gaussian distribution with mean $\log 3\log\log n$ and variance $(log 3)^2\log\log n$, which is an analog of the classical theorem of Erd\H os and Kac. And finally these results in all suggest that the behavior of $R(n;a)$ resembles the divisor function $d(n^2)$ in various aspects.    
\end{abstract}

\maketitle

\section{Introduction} \label{s1}
\noindent In the previous memoir of this series (see \cite{HV}) we studied the mean value
\begin{equation}
\label{e:eq3}
S(N;a)=\sum_{\substack{
n\leq N\\
(n,a)=1}
}R(n;a),
\end{equation}
of the number $R(n;a)$ of positive integer solutions to the Diophantine equation
\begin{equation} 
\label{e1}
\frac{a}{n}=\frac1x+\frac1y.
\end{equation}
Here we extend our investigation to the second moment and some consequences thereof. 

\begin{thm} 
\label{t1}
For fixed positive integer $a$, we have, for every $N\in\mathbb N$ with $N\ge 2$, 
$$\sum_{\substack{
n\le N\\
(n,a)=1}
}\Bigg|R(n;a) - \frac1{\phi(a)}\sum_{
\substack{\chi\bmod{a}\\
\chi^2=\chi_{{}_{\scriptstyle 0}}}
}\bar\chi(-n)\sum_{u|n^2}\chi(u)\Bigg|^2
\ll_a N\log^2 N,$$
where $\ll_a$ indicates that the implicit constant depends at most on $a$, and where $\chi_{{}_{\scriptstyle 0}}$ denotes the principal character modulo $a$. 
\end{thm}
In the character sum here the term $\chi=\chi_{{}_{\scriptstyle 0}}$ contributes an amount $d(n^2)$ where $d$ is the divisor function and we can expect that this is the dominant contribution on average.  Thus as a consequence of the Erd\H os--Kac theorem, just as for the divisor function $d(n)$, one can anticipate that $\log R(n;a)$ has a Gaussian distribution.  As a first approximation we establish the normal order of $\log R(n;a)$.
\begin{thm} \label{t2}
When $a$ is fixed, the normal order of $\log R(n;a)$ as a function of $n$ is $(\log 3)\log\log n$.
\end{thm}
\par
Let
$$\Phi(z):=\frac1{\sqrt{2\pi}}\int_{-\infty}^{z}e^{-\frac{t^2}{2}}dt.$$  Then with a little more work we can establish the full distribution.
\begin{thm} \label{t3}
For fixed positive integer $a$, we have
$$\lim_{N\to\infty}\frac1N\rm{card}\left\{n\le N:\frac{\log R(n;a)-(\log 3)\log\log n}{(\log 3)\sqrt{\log\log n}}\le z\right\}=\Phi(z).$$
\end{thm}

For completeness we also establish the mean square of $R(n;a)$ for fixed $a$. Since $R(n;a)$ resembles quite closely the divisor function $d(n^2)$ in many aspects, we expect that their mean squares share the same order of magnitude. Thus the following theorem can be compared with the asymptotic formula 
$$\sum_{n\le N}d^2(n^2)=N P_8(\log N)+O(N^{1-\delta})$$
which holds for some $\delta>0$ and with $P_8(\cdot)$ a polynomial of degree $8$.
\begin{thm} \label{t4}
Let $a$ be a fixed positive integer and $\varepsilon>0$. Then
$$\sum_{\substack{n\le N\\(n,a)=1}}R(n;a)^2=N P_8(\log N;a)+O_a(N^{35/54+\varepsilon})$$
where $P_8(\cdot\;;a)$ is a degree 8 polynomial with coefficients depending on $a$, and its leading coefficient is 
$$\frac1{8!a^2}\prod_{p|a}\left(1-\frac1p\right)^7\prod_{p\nmid a}\left(1-\frac6p+\frac1{p^2}\right)\left(1-\frac1p\right)^6.$$
\end{thm}

The error term in the theorem above is closely related to the generalised divisor problem, and in particular depends on a mean value estimate for the ninth moment of Dirichlet L-functions $L(s,\chi)$ modulo $a$ inside the critical strip. As is easily verified, the error can be improved to $O_a(N^{1/2+\varepsilon})$ under the assumption of the generalised Lindel\"{o}f Hypothesis.

\section{Proof of Theorem \ref{t1}} \label{s2}
\noindent We rewrite equation \eqref{e1} in the form
$$(ax-n)(ay-n)=n^2.$$
After the change of variables $u=ax-n$ and $v=ay-n$, it follows that $R(n;a)$ is the number of ordered pairs of natural numbers $u$, $v$ such that $uv=n^2$ and $u\equiv v\equiv -n\pmod{a}$.\par
Under the assumption that $(n,a)=1$, $R(n;a)$ can be reduced further to counting the number of divisors $u$ of $n^2$ with $u\equiv -n\pmod{a}$.  Now the residue class $u\equiv -n \pmod{a}$ is readily isolated {\it via} the orthogonality of the Dirichlet characters $\chi$ modulo $a$.  Thus we have
\begin{equation}
R(n;a) = \frac{1}{\phi(a)} \sum_{\substack{
\chi\\
\mod{a}
}} \Bar{\chi}(-n)\sum_{u|n^2}\chi(u),
\label{e2.5}
\end{equation}
where the condition $(n,a)=1$ is taken care of by the character $\bar{\chi}(n)$.\par
Hence
\begin{eqnarray*}
&&\sum_{\substack{n\le N\\
(a,n)=1
}}\Bigg|R(n;a)-\frac1{\phi(a)}\sum_{\substack{
\chi\,\mathrm{mod}\,\,a\\
\chi^2=\chi_{{}_{\scriptstyle 0}}
}}\bar\chi(-n)\sum_{u|n^2}\chi(u)\Bigg|^2\\
&\ll_a&\sum_{
\substack{n=1\\
(a,n)=1
}}^{\infty} e^{-n/N}\Bigg| \sum_{
\substack{\chi\,\mathrm{mod}\,\,a\\
\chi^2\not=\chi_{{}_{\scriptstyle 0}}
}}\bar\chi(-n)\sum_{u|n^2}\chi(u)\Bigg|^2\\
&=&\sum_{
\substack{\chi_{{}_{\scriptstyle 1}}\mathrm{mod}\,\,a\\
\chi_{{}_{\scriptstyle 1}}^2\not=\chi_{{}_{\scriptstyle 0}}
}}\,\, \sum_{
\substack{\chi_{{}_{\scriptstyle 2}}\mathrm{mod}\,\,a\\
\chi_{{}_{\scriptstyle 2}}^2\not=\chi_{{}_{\scriptstyle 0}}
}}  \bar\chi_{{}_{\scriptstyle 1}}\chi_{{}_{\scriptstyle 2}}(-1)\sum_{n=1}^{\infty} \sum_{u|n^2} \sum_{v|n^2} 
\chi_{{}_{\scriptstyle 1}}(u) \bar\chi_{{}_{\scriptstyle 2}}(v) \bar\chi_{{}_{\scriptstyle 1}}\chi_{{}_{\scriptstyle 2}}(n)e^{-n/N},
\end{eqnarray*}
where $\chi_{{}_{\scriptstyle 0}}$ denotes the principal character modulo $a$.
In order to evaluate the sum over $n$, we analyze the Dirichlet series
$$
f_{\chi_{{}_{\scriptstyle 1}},\chi_{{}_{\scriptstyle 2}}}(s):=\sum_{n=1}^{\infty} \sum_{u|n^2}\sum_{v|n^2}
\chi_{{}_{\scriptstyle 1}}(u)\bar\chi_{{}_{\scriptstyle 2}}(v)\bar\chi_{{}_{\scriptstyle 1}}\chi_{{}_{\scriptstyle 2}}(n) n^{-s}.
$$
The condition $u|n^2$ can be written as $u_1u_2^2|n^2$ with $u_1$ squarefree, i.e. $u_1u_2|n$, and likewise for $v|n^2$.  Thus 
\begin{equation}\label{e2.2}
f_{\chi_{{}_{\scriptstyle 1}},\chi_{{}_{\scriptstyle 2}}}(s) = \sum_{m=1}^{\infty}\frac{\bar\chi_{{}_{\scriptstyle 1}}\chi_{{}_{\scriptstyle 2}}(m)}{m^s}\sum_{d=1}^{\infty}\frac{F(d)}{d^s}
\end{equation}
where 
$$F(d)=\sum_{\substack{
u_1,u_2,v_1,v_2\\
d=[u_1u_2,v_1v_2]
}} \mu^2(u_1) \mu^2(v_1) \chi_{{}_{\scriptstyle 1}}(u_1u_2^2) \bar\chi_{{}_{\scriptstyle 2}}(v_1v_2^2) \bar\chi_{{}_{\scriptstyle 1}}\chi_{{}_{\scriptstyle 2}}(d).$$
The function $F$ is multiplicative and so the inner sum above is 
\begin{equation}\label{e2.3}
\prod_{p}\left(1+\sum_{k=1}^{\infty}F(p^k)p^{-ks}\right),
\end{equation}
where 
\begin{equation}\label{e2.4}
F(p^k)=\sum_{\substack{
u_1,u_2,v_1,v_2\\ [u_1u_2,v_1v_2]=p^k
}} \mu^2(u_1) \mu^2(v_1) \chi_{{}_{\scriptstyle 1}}(u_1u_2^2) \bar\chi_{{}_{\scriptstyle 2}}(v_1v_2^2) \bar\chi_{{}_{\scriptstyle 1}}\chi_{{}_{\scriptstyle 2}}(p^k).
\end{equation}
In particular we have
$$F(p)=\chi_{{}_{\scriptstyle 0}}(p)+\sum_{\chi\in \mathcal{X} \backslash\{\bar\chi_{{}_{\scriptstyle 1}}\chi_{{}_{\scriptstyle 2}}\}}\chi(p),$$
where $\mathcal{X} =\{\chi_{{}_{\scriptstyle 1}}, \chi_{{}_{\scriptstyle 2}}, \chi_{{}_{\scriptstyle 1}}\chi_{{}_{\scriptstyle 2}}, \chi_{{}_{\scriptstyle 1}}\bar\chi_{{}_{\scriptstyle 2}}, \bar\chi_{{}_{\scriptstyle 1}}, \bar\chi_{{}_{\scriptstyle 2}}, \bar\chi_{{}_{\scriptstyle 1}}\bar\chi_{{}_{\scriptstyle 2}}, \bar\chi_{{}_{\scriptstyle 1}}\chi_{{}_{\scriptstyle 2}}\}$ (and the entries are considered to be formally distinct), and
$$|F(p^k)|\le 8k.$$
Thus the Dirichlet series $f$ converges absolutely for $\sigma>1$ and  
\begin{equation}\label{e2.1}
f_{\chi_{{}_{\scriptstyle 1}},\chi_{{}_{\scriptstyle 2}}}(s)=G_{\chi_{{}_{\scriptstyle 1}},\chi_{{}_{\scriptstyle 2}}}(s)L(s,\chi_{{}_{\scriptstyle 0}}) \prod_{\chi\in \mathcal X}L(s,\chi),
\end{equation}
where $G_{\chi_{{}_{\scriptstyle 1}},\chi_{{}_{\scriptstyle 2}}}(s)$ is a function which is analytic in the region $\Re{s}>1/2$ and satisfies
$$G(s)\ll 1 \quad (\sigma\ge \textstyle\frac12+\delta)$$
for any fixed $\delta>0$.  As $\chi_{{}_{\scriptstyle 1}},\chi_{{}_{\scriptstyle 2}}$ are not characters of order 1 or 2, $f_{\chi_{{}_{\scriptstyle 1}},\chi_{{}_{\scriptstyle 2}}}(s)$ has a triple pole at $s=1$ when $\chi_{{}_{\scriptstyle 1}}=\chi_{{}_{\scriptstyle 2}}$ or $\chi_{{}_{\scriptstyle 1}}\chi_{{}_{\scriptstyle 2}} = \chi_{{}_{\scriptstyle 0}}$, and a simple pole otherwise.  By Corollary 1.17 and Lemma 10.15 of \cite{MV}, for fixed $a$, 
$$L(s,\chi)-\frac{E(\chi)\phi(a)}{a(s-1)}\ll 2+|t|$$
uniformly for $\sigma\ge \frac12$ where $E(\chi)$ is $1$ when $\chi=\chi_{{}_{\scriptstyle 0}}$ and $0$ otherwise.  Hence by (5.25) of \cite{MV}
$$\sum_{n=1}^{\infty} \sum_{u|n^2} \sum_{v|n^2} 
\chi_{{}_{\scriptstyle 1}}(u) \bar\chi_{{}_{\scriptstyle 2}}(v) \bar\chi_{{}_{\scriptstyle 1}}\chi_{{}_{\scriptstyle 2}}(n)e^{-n/N} = \frac1{2\pi i}\int_{\theta-i\infty}^{\theta+i\infty}f_{\chi_{{}_{\scriptstyle 1}},\chi_{{}_{\scriptstyle 2}}}(s) N^s \Gamma(s) ds,$$
where $\theta>1$.
Since the gamma function decays exponentially fast on any vertical line we may move the vertical path to the $\frac34$--line picking up the residue of the integrand at $s=1$.  The  residue contributes an amount
$$\ll N(\log N)^2$$
and the new vertical path contributes
$$\ll N^{\frac34}.$$ 
This completes the proof of Theorem \ref{t1}.

\section{Proof of Theorem \ref{t2}} \label{s3}
\noindent By Theorem \ref{t1}, we expect that for almost all $n$ with $(a,n)=1$, $R(n;a)$ is close to $$\frac1{\phi(a)}\sum_{\substack{\chi\,\mathrm{mod}\,\,a\\\chi^2=\chi_{{}_{\scriptstyle 0}}}}\bar\chi(-n)\sum_{u|n^2}\chi(u).$$ 
Thus we need to examine the contribution from the characters modulo $a$ of order 1 and 2.  For general $a$, there may be many quadratic characters modulo $a$.  Nevertheless we believe that the major contribution to the sum above comes from
the principle character, and this is of size 
$$\frac{d(n^2)}{\phi(a)}.$$  
Thus, for fixed $a$, $\log R(n;a)$ should have the normal order of $\log d(n^2)$, namely $(\log 3)\log\log n$.  When $(n,a)>1$ we have 
\begin{equation}
R(n;a)=R(n/(n,a);a/(n,a))
\label{e2.6}
\end{equation}
and so we can expect that the general case follows from the special case $(n,a)=1$.

Before embarking on the proof of Theorem \ref{t2}, we state a lemma.
We define, for any quadratic character $\chi$,
$$\Omega_\chi(n)=\mathrm{card}\left\{p,k: k\ge 1, p^k|n, \chi(p^k)=1\right\}.$$

\begin{lem} \label{l1}
Suppose that $\chi$ is a quadratic character to a fixed modulus $a$ and that $N\ge 3$.  Then
$$\sum_{n\le N}\Bigg(\Omega_\chi(n)-\frac{1}{2}\log\log N\Bigg)^2\ll N\log\log N$$
and
$$\sum_{1<n\le N}\Bigg(\Omega_\chi(n)-\frac{1}{2}\log\log n\Bigg)^2\ll N\log\log N.$$
\end{lem}

\begin{proof}
The proof follows in the same way as Tur\'an's theorem (see Theorem 2.12 of \cite{MV}) on observing that 
$$\sum_{\substack{
p\le N\\
\chi(p)=1}} \frac1p =\frac12\log\log N +O(1)$$
and this is readily deduced from Corollary 11.18 of \cite{MV}.
\end{proof}
\par
It is an immediate consequence of the above lemma that $\Omega_\chi(n)$ has normal order $\frac12\log\log n$.  In particular, for any fixed $\varepsilon>0$, for almost all $n$,
$$3^{\Omega_{\chi}(n)}< 3^{(\frac12+\varepsilon)\log\log n}.$$
Now, for any quadratic character $\chi$ modulo $a$, let
$$g_{\chi}(n)=\sum_{u|n^2}\chi(u).$$
This is
$$\prod_{p^k\parallel n}\left(1+\chi(p)+\chi^2(p)+\cdots+\chi^{2k}(p)\right).$$
When $\chi(p)=-1$ the general factor is $1$, and when $\chi(p)=1$ it is $2k+1$.  Hence
$$0<g_{\chi}(n)\le 3^{\Omega_{\chi}(n)}.$$
Thus for any fixed $\varepsilon>0$, for every quadratic character modulo $a$, for almost all $n$, 
\begin{equation}
g_{\chi}(n) < (\log n)^{(\frac12\log 3+\varepsilon)}.
\label{e2.7}
\end{equation}
\par
Let $$r(n;a)=\frac1{\phi(a/(n,a))}\sum_{\substack{\chi\,\mathrm{mod}\,\,a/(n,a)\\\chi^2=\chi_{{}_{\scriptstyle 0}}}}\bar\chi(-n/(n,a))g_\chi(n/(n,a)).$$
Since $R(n;a)=R(n/(n,a);a/(n,a))$, it follows by Theorem \ref{t1} that
$$\sum_{n\le N}(R(n;a)-r(n;a))^2=\sum_{d|a}\sum_{\substack{
m\le N/d\\
(m,a/d)=1
}} (R(m;a/d)-r(m;a/d))^2
\ll N(\log N)^2.$$
Hence, for any fixed $\varepsilon>0$, for almost all $n$ we have 
$$|R(n;a)-r(n;a)|<(\log n)^{1+\varepsilon}.$$
Therefore, by \eqref{e2.7}, for almost all $n$, 
\begin{equation}
\left|
R(n;a)-\frac{d\left(
(n/(a,n))^2
\right)}{\phi(a/(a,n))}
\right|< (\log n)^{1+2\varepsilon}.
\label{e2.8}
\end{equation}
Now $3\le d(p^{2k})=1+2k\le 3^k$.  Hence 
\begin{equation}
3^{\omega(n)-\omega(a)}\le d\left(
(n/(a,n))^2
\right) \le 3^{\Omega(n)}
\label{e2.9}
\end{equation}
and it follows that
$$(\log n)^{\log 3-\varepsilon} < \frac{d\left(
(n/(a,n))^2
\right)}{\phi(a/(a,n))} < (\log n)^{\log 3+\varepsilon}$$
for almost all $n$.  Theorem \ref{t2} now follows.

\section{Proof of Theorem \ref{t3}} \label{s4}
\noindent By (\ref{e2.8}) and (\ref{e2.9}), for every fixed $\varepsilon>0$, for almost all $n$,
$$\frac{3^{\omega(n)}}{\phi(a/(a,n))} - (\log n)^{1+\varepsilon} < R(n;a) < 3^{\Omega(n)} +(\log n)^{1+\varepsilon}.$$
Moreover, for almost all $n$ we have $\Omega(n)\ge \omega(n)>(1-\varepsilon)\log\log n$.  Hence for any $\delta$ with $0<\delta< \log 3-1$ we have, for almost all $n$
$$3^{\omega(n)-\omega(a)- \log\phi(a/(a,n)) }\exp(-(\log n)^{-\delta}) < R(n;a) < 3^{\Omega(n)}\exp((\log n)^{-\delta})$$
and so
$$3^{\omega(n)}\exp(-\varepsilon\sqrt{\log\log n})<R(n;a)< 3^{\Omega(n)}\exp(\varepsilon\sqrt{\log\log n}).$$
Let
$$S(N;z)=\mathrm{card}\left\{n\le N:
\frac{\log R(n;a) - (\log 3)\log\log n}{\log 3\sqrt{\log\log n}} \le z
\right\},$$
$$S_-(N;z)= \mathrm{card}\left\{n\le N:
\frac{\Omega(n) - \log\log n}{\sqrt{\log\log n}} \le z
\right\}$$
and
$$S_+(N;z)= \mathrm{card}\left\{n\le N:
\frac{\omega(n) - \log\log n}{\sqrt{\log\log n}} \le z
\right\},$$
Then for a non-negative monotonic function $\eta(n)$ tending to $0$ sufficiently slowly as $N\rightarrow\infty$ we have 
$$-\eta(N)N +S_-(N;z-\varepsilon)< S(N;z)< \eta(N)N + S_+(N;z+\varepsilon).$$
Hence, by the Erd\H os--Kac theorem (see, for example Theorem 7.21 and Exercise 7.4.4 of \cite{MV}), 
$$\Phi(z-\varepsilon) \le \liminf_{N\rightarrow\infty} N^{-1} S(N;z) \le \limsup_{N\rightarrow\infty}N^{-1} S(N;z) \le \Phi(z+\varepsilon).$$
The theorem now follows from the continuity of $\Phi$.

\section{Proof of Theorem \ref{t4}} \label{s5}
\noindent By a similar discussion to that in \S\ref{s2}, we can show that the generating Dirichlet series for $R(n;a)^2$ is
$$\sum_{\substack{n=1\\(n,a)=1}}^{\infty}\frac{R(n;a)^2}{n^s}=\frac1{\phi(a)^2}\sum_{\substack{\chi_{{}_{\scriptstyle 1}},\chi_{{}_{\scriptstyle 2}}\\\mathrm{mod}\,\, a}}\bar\chi_{{}_{\scriptstyle 1}}\chi_{{}_{\scriptstyle 2}}(-1)
f_{\chi_{{}_{\scriptstyle 1}},\chi_{{}_{\scriptstyle 2}}}(s),$$ 
where $f_{\chi_{{}_{\scriptstyle 1}},\chi_{{}_{\scriptstyle 2}}}(s)$ is analytic in the region $\Re s>1/2$ and is given by \eqref{e2.1}.  Here $f_{\chi_{{}_{\scriptstyle 1}},\chi_{{}_{\scriptstyle 2}}}(s)$ has a pole at 1 of order at least 1, and as high as 9 exactly when $\chi_{{}_{\scriptstyle 1}}$ and $\chi_{{}_{\scriptstyle 2}}$ are equal to the principle character $\chi_{{}_{\scriptstyle 0}}$. Now on applying Perron's formula, we have for $\theta=1=1+1/\log(2N)$,
\begin{equation}
\sum_{\substack{n\le N\\(n,a)=1}}R(n;a)^2=\frac1{\phi(a)^2}\sum_{\substack{
\chi_{{}_{\scriptstyle 1}},\chi_{{}_{\scriptstyle 2}}\\\bmod a}} \frac{\bar\chi_{{}_{\scriptstyle 1}}\chi_{{}_{\scriptstyle 2}}(-1)}{2\pi i}\int_{\theta-i T}^{\theta+i T}f_{\chi_{{}_{\scriptstyle 1}},\chi_{{}_{\scriptstyle 2}}}(s)\frac{N^s}s d s+O_a(N^{1+\varepsilon}/T).
\label{e2.10}
\end{equation}

Since we are shooting for the asymptotics for the the mean square, smoothing factors of the kind used in section \ref{s2} are best avoided.  Since the integrand includes a product of nine $L$--functions, we cannot expect to to be able to move the vertical integral path too close to the 1/2-line in the current state of knowledge. Nevertheless, the following result of Meurman \cite{Me1} which extends Heath-Brown's theorem \cite{H-B} on the twelfth power moment of the Riemann zeta function to Dirichlet $L$--functions, provides a starting point for the analysis.

\begin{lem}\label{l6} 
\[
\sum_{\chi \bmod a}\int_{-T}^{T}|L(\textstyle{\frac12}+i t,\chi)|^{12}dt\ll a^3 T^{2+\varepsilon},
\]
where $\varepsilon>0$, $a\ge1$ and $T\ge 2$.
\end{lem}

Then, adapting the argument of Chapter 8 of Ivi\'{c} \cite{Iv} for the Riemann zeta function to Dirichlet L-functions establishes the following.  

\begin{lem}\label{l7}
\[
\int_{-T}^{T}|L({\textstyle\frac{35}{54}}+i t,\chi)|^{9}dt\ll_a T^{1+\varepsilon},
\]
where $\varepsilon>0$, $\chi$ is a fixed Dirichlet character modulo $a\ge1$ and $T\ge 2$.
\end{lem} 
If one utilizes the sharpest estimates for the underlying exponential sums, Lemma \ref{l7} is susceptible to slight improvements.

Now, we move the vertical integral path in (\ref{e2.10}) to the $35/54$-line, picking up the residue of the integrand at $1$.  Thus
\begin{align*}
\int_{\theta-i T}^{\theta+i T}f_{\chi_{{}_{\scriptstyle 1}},\chi_{{}_{\scriptstyle 2}}}(s)\frac{N^s}s d s
=&\int_{\theta-i T}^{35/54-i T}f_{\chi_{{}_{\scriptstyle 1}},\chi_{{}_{\scriptstyle 2}}}(s)\frac{N^s}s d s
+\int_{35/54+i T}^{\theta+i T}f_{\chi_{{}_{\scriptstyle 1}},\chi_{{}_{\scriptstyle 2}}}(s)\frac{N^s}s d s\\
&+\int_{35/54-i T}^{35/54+i T}f_{\chi_{{}_{\scriptstyle 1}},\chi_{{}_{\scriptstyle 2}}}(s)\frac{N^s}s d s+\text{Res}_{s=1}\left(f_{\chi_{{}_{\scriptstyle 1}},\chi_{{}_{\scriptstyle 2}}}(s)\frac{N^s}s\right)
\end{align*}

Here, in order to deal with the contribution from the horizontal integrals, we cannot afford to use the crude convexity bounds on Dirichlet L-functions, due to the large number of $L$--functions in the integrand.  Fortunately, a sharper bound  has been established by Pan \& Pan in Theorem 24.2.1 of \cite{PP}. 

\begin{lem}\label{l8}
Let $l\ge3$, $L=2^{l-1}$ and $\sigma_l=1-l(2L-2)^{-1}$. Then when $\sigma\ge\sigma_l$
$$L(\sigma+it,\chi)\ll_a |t|^{1/(2L-2)}\log|t|$$
holds uniformly for $|t|\ge2$.
\end{lem}
When $l=3$ we obtain
$$L(\sigma+it,\chi)\ll_a |t|^{1/6}\log|t|$$
uniformly for $|t|\ge 2$ and $\sigma\ge \frac12$, and when $l=4$,
$$L(\sigma+it,\chi)\ll_a |t|^{1/{14}}\log|t|$$
uniformly for $\sigma \ge 5/7$.  Thus, by the convexity principle for Dirichlet series, 
$$L(\sigma+it,\chi)\ll_a |t|^{\mu(\sigma)+\varepsilon}$$
uniformly for $|t|\ge 2$ and $\sigma\ge\frac12$ where
$$\mu(\sigma)=\begin{cases}\textstyle\frac16-\frac49(\sigma-\frac12)&\text{ when }\frac12\le\sigma\le\frac57,\\
\textstyle\frac{1-\sigma}4&\text{ when }\frac57<\sigma\le 1,\\
0&\text{ when }1<\sigma.
\end{cases}$$
We note that $\mu(\frac{35}{54})=\frac{49}{486}<\frac19$ and $\mu(\frac57)=\frac1{14}$.

Now the horizontal paths contribute
$$\ll \int_{35/54}^{1+\varepsilon}N^{\sigma}|f_{\chi_{{}_{\scriptstyle 1}},\chi_{{}_{\scriptstyle 2}}}(\sigma+i T)|T^{-1}d \sigma.$$
and this is
$$\ll\max_{35/54\le\sigma\le1+\varepsilon}N^{\sigma}T^{9\mu(\sigma)-1+
\varepsilon},$$
and by the piecewise linearity of $\sigma$ and $\mu(\sigma)$ this is
$$\ll N^{1+\varepsilon}T^{-1} + N^{5/7}T^{9\mu(5/7)-1+\varepsilon} + N^{35/54}T^{9\mu(35/54)-1+\varepsilon}.$$  
When $T=N$ this is 
$$\ll N^{35/54+\varepsilon}.$$
On the other hand, by Lemma \ref{l7} the vertical path also contributes
$$\ll N^{35/54+\varepsilon}.$$  

The main term comes from the residual contributions, which, in the case that $\chi_{{}_{\scriptstyle 1}}=\chi_{{}_{\scriptstyle 2}}=\chi_{{}_{\scriptstyle 0}}$, is
$N P_8(\log N;a)$ where $P_8(\cdot\;;a)$ is a polynomial of degree 8 whose coefficients depend on $a$. Notice that for other choices of $\chi_{{}_{\scriptstyle 1}}$ and $\chi_{{}_{\scriptstyle 2}}$, the residual contribution gives a polynomial of $\log N$ of lower degree than above.

For the leading coefficient, we need more precise information about $f_{\chi_{{}_{\scriptstyle 0}},\chi_{{}_{\scriptstyle 0}}}$.  By \eqref{e2.2}, \eqref{e2.3} and \eqref{e2.4} we have 
$$f_{\chi_{{}_{\scriptstyle 0}},\chi_{{}_{\scriptstyle 0}}}=L(s,\chi_{{}_{\scriptstyle 0}})\prod_{p\nmid a}\left(1+\sum_{k=1}^{\infty}\frac{8k}{p^{k s}}\right)=L(s,\chi_{{}_{\scriptstyle 0}})^9 \prod_{p\nmid a}(1+6p^{-s}+p^{-2s})(1-p^{-s})^6,$$
from which the leading coefficient is readily deduced. This completes the proof of Theorem \ref{t4}.
\bigskip
\par
In conclusion we remark that a concomitant argument will give
$$\sum_{\substack{n\le N\\(n,a)=1}}R(n;a)^k=N P_{3^k-1}(\log N;a)+O_a\left(N^{\alpha_k+\varepsilon}\right)$$
for any $\varepsilon>0$, where $P_{3^k-1}(\cdot\;;a)$ is a polynomial of degree $3^k-1$ whose coefficients depend on $a$ and $\alpha_k$ is a constant that depends on the best $3^k$-th power moment estimates for $L(s,\chi)$ in the critical strip and the quantity $\mu(\sigma)$ defined above. This question is closely related to the generalised divisor problem, and one is referred to Chapter 13 in Ivi\'{c} \cite{Iv} for more details.


\begin{thebibliography}{xxxx}
\bibitem[H-B]{H-B} D. R. Heath-Brown, \emph{Mean values of the zeta-function and divisor problems}, Recent Progress in Analytic Number Theory, Vol. 1 (Durham, 1979), Acad. Press, London-New York, 1981, pp. 115-119. 
\bibitem[HV]{HV} J. J. Huang and R. C. Vaughan, \emph{Mean value theorems for binary Egyptian fractions}, J. Number Theory 131 (2011), no. 11, 1641-1656.
\bibitem[Iv]{Iv} A. Ivi\'{c}, \emph{The Theory of the Riemann Zeta-Function with Applications}. A Wiley-Interscience Publication. John Wiley \& Sons, Inc., New York, 1985. xvi+517 pp. 
\bibitem[Me1]{Me1} T. Meurman, \emph{The mean twelfth power of Dirichlet $L$-functions on the critical line}, Ann. Acad. Sci. Fenn. Ser. A I Math. Dissertationes No. 52 (1984), 44 pp.
\bibitem[Me2]{Me2} T. Meurman, \emph{A generalization of Atkinson's formula to $L$-functions}, Acta Arith. 47 (1986), no. 4, 351--370. 
\bibitem[MV]{MV} H. L. Montgomery and R. C. Vaughan, \emph{Multiplicative Number Theory I. Classical Theory}, Cambridge University Press, 2007.
\bibitem [PP]{PP} C. D. Pan and C. B. Pan, \emph{Foundation to Analytic Number Theory (Chinese)}, Science Press, Beijing, 1991. 


\end{thebibliography}
\end{document}